# Polymultisets, Multisuccessors, and Multidimensional Peano Arithmetics

**Alexander Chunikhin**
Palladin Institute of Biochemistry
National Academy of Sciences of Ukraine
E-mail address: alchun@ukr.net

**Abstract.** The goal of this paper is introduction of a concept of natural multidimensional numbers and to construct a generalized Peano arithmetic of these multidimensional numbers. For this purpose we define a polymultiset as a special set-like form of m-ary multirelation. In addition, multisuccessor and multipredecessor functions on polymultisets are determined. By introducing Peano-like axioms for natural multidimensional numbers, we build a commutative semiring of 2-numbers $\langle \mathbf{P_2}, +, \cdot, \mathbf{0_2}, 1_{00} \rangle$.

## 1. Introduction

Multisets and multi-relations play an important role in Database and Data Mining [2, 7].

**Definition 1.1.** [1, 5] Let A and B be sets. A *multiset relation* (*multi-relation*) on A and B is a multiset M: $A \times B \to N_0$, $N_0 = \{0, 1, 2, \ldots\}$, where M(a, b) defines the multiplicity of the pair $(a, b) \in A \times B$.

Let $\mathbf{G} = \{G_1, G_2, \ldots, G_m\}$ be a set of domains that characterize particular properties (attributes) of some system (a complex object). Suppose that each domain $G_i$ consists of a finite number of distinct elements $G_i = \{g_{1i}, g_{2i}, \ldots, g_{ki}\}$.

Let us take the Cartesian product $\mathbf{G} = G_1 \times G_2 \times \ldots \times G_m$ of domains, which is called the domain base.

**Definition 1.2.** The *m-ary multi-relation* $R_m$ is defined as $R_m: G_1 \times G_2 \times \ldots \times G_m \to N_0$, where $G_i$ is the domain of *i-th* attribute.

Elements from the domain base $\mathbf{G}$ have the form $I_{kl\ldots r} = (g_{1k}, g_{2l}, \ldots, g_{mr})$ and are called *polyments*.

Denote the occurrence (degree of multiplicity) of the polymlement $I_{kl\ldots r}$ by $\alpha_{k\ldots r}$ $= \# (g_{1k}, \ldots, g_{mr}, R_m)$, $\alpha_{k\ldots r} \in N_0$. Then a collection of polymlements with degrees of multiplicity corresponding to them forms a (multi)set-similar structure. It is natural to suggest the following logical chain: m-ary multirelation – special set-similar structure – the number of a special form.

## 2. The PolyXsets concept

The first step of the logical chain mentioned above is to define a multi-relation as a special set-similar structure.

**Definition 2.1.** Given a set of domains $\mathbf{G} = \{G_1, G_2, \ldots, G_m\}$, let us define the m-dimensional *polyXset* $A_m$ as a pair $(\mathbf{I}, \alpha)$ where $\mathbf{I}$ is a set of m-dimensional polyments $\{I_{k\ldots r}\}$ and $\alpha\colon \mathbf{I} \to \mathfrak{I}$ is a function from $\mathbf{I}$ to a numerical structure $\mathfrak{I}$. The number $\alpha_{k\ldots r} = \alpha(I_{k\ldots r})$ is called a degree of multiplicity (occurrence) of the polyment $I_{k\ldots r} \in \mathbf{I}$ in $A_m$:

$$A_m = (\mathbf{I}, \alpha) = <(I_{k\ldots r}, \alpha_{k\ldots r})>_{k\ldots r}.$$

The polyment $I_{kl\ldots r}$ together with the degree of multiplicity $\alpha_{k\ldots r}$ is called a *component* of the polyXset $A_m$.

If $\mathfrak{I} = N_0$, then the polyXset is called *polymultiset* or *polymset*.
If $\mathfrak{I} = [0, 1]$, then the polyXset is called *polyfuzzyset* or *polyfset*.
If $\mathfrak{I} = Z$, then the polyXset is called *polyhybridset* or *polyhset*.

Below, we exactly describe polymsets.

**Example 2.2.** Let us consider the polymset $A_3$, which has three-component basis of properties $\mathbf{G} = \{G_1, G_2, G_3\}$:
A body form $G_1 = \{\text{cube, pyramid, sphere, cone, cylinder}\}$;
A material $G_2 = \{\text{metal, plastic, paper}\}$;
A color $G_3 = \{\text{black, white, red, yellow, green, blue}\}$.
Suppose the polymset $A_3$ consists of one blue plastic cube, three black metal spheres, seven yellow paper cones, and five black metal cylinders.
For elements of each property $G_j$, we will enter an arbitrary numbering, i.e. it is comparable to the first index of each element $g_{ij} \in G_j$ natural number, for example, the element's number in the list.
Then "a blue plastic cube" is represented in the basis $\mathbf{G}$ as $I_{126} = (g_{11}, g_{22}, g_{63})$, "black metal sphere" is represented as $I_{311} = (g_{31}, g_{12}, g_{13})$.
Besides, $I_{126}$ ("a blue plastic cube") is the polyment of the polymset $A_3$ in the basis $\mathbf{G}$, and $(I_{434}, 7)$ ("seven yellow paper cones") is the component of the polymset $A_3$, and $\alpha_{434} = 7$.
Thus, the polymset $A_3$ can be represented in the form

$$A_3 = <(I_{126}, 1), (I_{331}, 3), (I_{434}, 7), (I_{511}, 5)>.$$

**Definition 2.3.** The *support* of polymset $A_m$ (Supp $A_m$) is defined as a polymset of individual copies of all polyments from $A_m$:

$$\text{Supp } A_m = \{I_{k\ldots r} \mid \chi_{k\ldots r}(A_m) = 1\},$$

where $\chi_{k\ldots r}(A_m)$ is a polymset characteristic function such that $\chi_{k\ldots r}(A_m) = 1$ if $\alpha_{k\ldots r} \neq 0$ and $\chi_{k\ldots r}(A_m) = 0$ if $\alpha_{k\ldots r} = 0$.

**Definition 2.4.** The *cardinality* of polymset $A_m$ is defined as the total number of all its polyments copies: Card $A_m = \sum_{k\ldots r} \alpha_{k\ldots r}$.

**Definition 2.5.** An *empty polymset* $\emptyset$ is defined as a polymset that does not have polyments.

$$\emptyset_m: \alpha_{k\ldots r} = 0, \forall k,\ldots, r.$$

**Definition 2.6.** The *height of polymset* $A_m$ is defined as the greatest degree of multiplicity for all polyment copies: $\text{hgt}(A_m) = \max_{k\ldots r} \alpha_{k\ldots r}$.

**Definition 2.7.** a) Polymset $A_m$ is called *constant* if for any $k, \ldots, r$, $\alpha_{k\ldots r} = \text{const}$.

b) A family of polymsets is called *n-bounded* if $\alpha_{ij\ldots k} \leq n$ for each polyment $I_{ij\ldots k}$ of polymset $A_m$, i.e., any polyment cannot occur more than $n$ times in $A_m$.

c) In general, it is possible to have such a restriction for each polyment: $\alpha_{ij\ldots k} \leq n_{ij\ldots k}$. A family of such polymsets is called *individually bounded.*

**Definition 2.8.** Two polymsets $A_m$ and $B_m$ are *equal* ($A_m = B_m$) if $I_{k\ldots r}(A_m) \equiv I_{k\ldots r}(B_m)$ and $\alpha_{k\ldots r} = \beta_{k\ldots r}, \forall k,\ldots, r$, i.e., the polymsets consist of identical polyments with equal degree of multiplicity. If the second condition does not take place, then the polymsets are called *similar*.

**Definition 2.9.** A polymset $B_m$ is a *subpolymset* of polymset $A_m$ ($B_m \subseteq A_m$) if $\beta_{k\ldots r} \leq \alpha_{k\ldots r}, \forall k,\ldots, r$. Also Card $B_m \leq$ Card $A_m$; Supp $B_m \subseteq$ Supp $A_m$.

**Definition 2.10.** Polymsets $A_m$ and $B_m$ are called:

(i) *equicardinal* if Card $A_m =$ Card $B_m$. Equal polymsets are equicardinal, but not vice versa.

(ii) *equidimensional* if dim $A_m =$ dim $B_m$. Equal polymsets are equidimensional, but not vice versa.

(iii) *equivalent* if they are equicardinal and equidimensional.

Hence a *set* is one-dimensional one-bounded polymset, and a *multiset* is one-dimensional polymset with natural multiplicities of polyments, and a *fuzzy set* is one-dimensional polymset with multiplicities from [0; 1].

### 3. Operations on polymsets

Operations mentioned below are defined for polymsets with equal domains only. Let $A_m = \langle(\alpha_{k\ldots r}, I_{k\ldots r})\rangle$, $B_m = \langle(\beta_{k\ldots r}, J_{k\ldots r})\rangle$, $\alpha_{k\ldots r}, \beta_{k\ldots r} \in N_0$ be polymsets.

**Definition 3.1.** The *union* of two polymsets $A_m$ and $B_m$ is a polymset $C_m$ denoted by $C_m = A_m \cup B_m$ such that $C_m = \langle(\gamma_{k\ldots r}, \Gamma_{k\ldots r} \mid \gamma_{k\ldots r} = \max(\alpha_{k\ldots r}, \beta_{k\ldots r}))\rangle$.

**Definition 3.2.** The *intersection* of two polymsets $A_m$ and $B_m$ is a polymset $D_m$ denoted by $D_m = A_m \cap B_m$ such that $D_m = \langle(\delta_{k\ldots r} \cdot \Gamma_{k\ldots r} \mid \delta_{k\ldots r} = \min(\alpha_{k\ldots r}, \beta_{k\ldots r}))\rangle$.

**Definition 3.3.** *Arithmetical addition* of two polymsets $A_m$ and $B_m$ is defined as a polymset $C_m$ denoted by $C_m = A_m + B_m$ such that $C_m = <(\gamma_{k...r} \cdot \Gamma_{k...r} \mid \gamma_{k...r} = \alpha_{k...r} + \beta_{k...r})>$.

**Definition 3.4.** *Arithmetical subtraction* of two polymsets $A_m$ and $B_m$ is defined as a polymset $C_m$ denoted by $C_m = A_m - B_m$ such that $C_m = <(\gamma_{k...r} \cdot \Gamma_{k...r} \mid \gamma_{k...r} = \max((\alpha_{k...r} - \beta_{k...r}), 0))>$.

**Definition 3.5.** *Symmetric difference* of two polymsets $A_m$ and $B_m$ is a polymset $C_m$ denoted by $C_m = A_m \Delta B_m$ such that $C_m = <(\gamma_{k...r} \cdot \Gamma_{k...r} \mid \gamma_{k...r} = |\alpha_{k...r} - \beta_{k...r}|)>$. At the same time, $\gamma_{k...r} = \max(\alpha_{k...r}, \beta_{k...r}) - \min(\alpha_{k...r}, \beta_{k...r})$. Moreover, $A_m \Delta B_m = (A_m - B_m) + (B_m - A_m)$.

**Definition 3.6.** The *reduction* (decrease of dimension) of polymset $A_m$ on domain $G_i$ is defined as the polymset $A_{m-1}$ such that $\dim A_{m-1} = m-1$, and multiplicities of polyments $I_{i...k}$ in the reduced polymset are defined as $\alpha_{i...k} = \sum_{j=1}^{q} \alpha_{ij...k}$, where $q =$ card $G_i$. The full reduction of the polymset $A_m$ (convolution on all *m* domains) is equivalent to getting of polymset cardinality: $\sum_{i} ... \sum_{l} \alpha_{i...l} = $ Card $A_m$.

**Definition 3.7.** The *production* (increase of dimension) of the polymset $A_m$ on domain $G_t$ is defined as polymset $A_{m+1}$ such that $\dim A_{m+1} = m+1$, and Card $A_{m+1} = $ Card $A_m$, and

$$\sum_{t=1}^{q} \alpha_{i...kt} = \alpha_{i...k}.$$

## 4. A principle of polymsets generation

The inductive principle of number generation [4, 5] plays a basic role in the definition of natural numbers.

Let us use this principle. We will show that an arbitrary polymset $A_m$ is formed (generated) from the empty polymset $\emptyset$ by using one or several special multisuccessor functions $Sc_{k...r}(.)$.

Thus, under the successing we will understand an addition to initial polymset one copy of the corresponding polyment:

$Sc_{k...r}(\emptyset) = <(I_{k...r}, 1)>$,
$Sc_{k...r}(<(I_{k...r}, 1)>) = Sc_{k...r}(Sc_{k...r}(\emptyset)) = Sc^2_{k...r}(\emptyset) = <(I_{k...r}, 2)>$,
...............................................
$Sc_{k...r}(Sc^{(\alpha-1)}_{k...r}(\emptyset)) = Sc_{k...r}(Sc_{k...r}...Sc_{k...r}(\emptyset)...) = Sc^{\alpha}_{k...r}(\emptyset) = <(I_{k...r}, \alpha)>$.

At the same time, application of the successor function to index combination that is different from the initial combination adds a new polyment to initial polymset.

$Sc_{p...q}(A_m) = A_m + <(I_{p...q}, 1)>$,
$Sc^{\alpha}_{p...q}(A_m) = A_m + <(I_{p...q}, \alpha)>$.

A process of arbitrary polymset generation can be described as a sequence:

$\emptyset$, $Sc_{i\ldots j}(\emptyset)$, $Sc_{p\ldots q}(Sc_{i\ldots j}(\emptyset))$, $Sc_{s\ldots t}(Sc_{p\ldots q}(Sc_{i\ldots j}(\emptyset)))$, … .

**Example 4.1.** For $\mathbf{G} = \{G_1, G_2, G_3\}$ (see example 2.2) we can obtain

$Sc_{126}(\emptyset) = \langle(I_{126}, 1)\rangle$ ("a blue plastic cube"),
$Sc_{126}(\langle(I_{126}, 1)\rangle) = Sc^2_{126}(\emptyset) = \langle(I_{126}, 2)\rangle$ ("two blue plastic cubes").

Meanwhile,

$Sc_{434}(\langle(I_{126}, 1)\rangle) = Sc_{434}(Sc_{126}(\emptyset)) = \langle(I_{126}, 1), (I_{434}, 1)\rangle$ ("one blue plastic cube" and (plus) "one yellow paper cone").

**Definition 4.2.** A polymset $B_m$ is called the *immediate successor* of the polymset $A_m$ iff $A_m$ is a subpolymset of $B_m$ and a cardinality of $B_m$ is more by one than a cardinality of $A_m$:

$B_m = Sc(A_m) \Leftrightarrow (A_m \subseteq B_m, \text{Card } B_m = \text{Card } A_m + 1)$.

This is a general definition of a succession. It does not consider a constructive component, i.e., what polyment from a set of possible ones has been added?
We have to specify the definition of successions.

**Definition 4.3.** A polymset $B_m$ is called the *immediate i…j-successor* (multisuccessor) of the polymset $A_m$ on the set of indexes $(i, \ldots, j)$ iff one copy of polyment $I_{i\ldots j}$ is added to the polymset $A_m$ (see def. 3.3): $B_m = Sc_{i\ldots j}(A_m) \Leftrightarrow B_m = A_m + I_{i\ldots j}$.

The multiplicity degree $\beta_{i\ldots j}$ of polyment $I_{i\ldots j}$ in polymset-successor $B_m$ increases by one:
$\beta_{i\ldots j} = \alpha_{i\ldots j} + 1$.

**Definition 4.4.** A polymset $B_m$ is called the *immediate predecessor* of polymset $A_m$ iff $B_m$ is a subpolymset of $A_m$ and a cardinality of $B_m$ is less by one than a cardinality of $A_m$:

$B_m = Pd(A_m) \Leftrightarrow (B_m \subseteq A_m, \text{Card } B_m = \text{Card } A_m - 1)$.

This is a general definition of the predecessor too.

**Definition 4.5.** A polymset $B_m$ is called *the immediate i…j-predecessor* (multipredecessor) of polymset $A_m$ on the set of indexes $(i, \ldots, j)$ iff it can be obtained from the polymset $A_m$ by removal of one copy of polyment $I_{i\ldots j}$ (see def.3.4): $B_m = Pd_{i\ldots j}(A_m) \Leftrightarrow B_m = A_m - I_{i\ldots j}$.

The multiplicity degree $\beta_{i\ldots j}$ of polyment $I_{i\ldots j}$ in polymset-predecessor $B_m$ decreases by one: $\beta_{i\ldots j} = \alpha_{i\ldots j} - 1$.
Hereinafter the multisuccessor $\equiv$ i…j-successor.

Note that for determination of arbitrary polymset over domain base $\mathbf{G} = \{G_1, G_2, \ldots, G_m\}$ it is sufficient to specify the degrees of all multiplicities $\alpha_{i\ldots j}$ ($\forall i, \ldots, j$) only. If we will enter a numbering on the set of elements of properties, then each polyment

$I_{k \ldots r}$ can be unequivocally define by the set of indexes $k, \ldots, r$, and the polymset form $A_m = <(I_{k \ldots r}, \alpha_{k \ldots r})>_{k \ldots r}$ (see def.2.1) can be reduced to the form $A_m = <\alpha_{k \ldots r}>_{k \ldots r}$.

So instead of operating with polymsets actually we can operate with the degrees of multiplicity $\alpha_{i \ldots j}$. For example 2.2 we obtain

$A_3 = <(I_{126}, 1), (I_{331}, 3), (I_{434}, 7), (I_{511}, 5)> = <1_{126}, 3_{331}, 7_{434}, 5_{511}>$.

It is useful to represent the multiindex components of polymset in multidimensional matrix form.

Agree not formally that a multidimensional matrix $\alpha_{i \ldots j} \rfloor$, $\forall \alpha_{i \ldots j} \in N_0$ be natural multidimensional number. In case of two dimensions (domain base $G = \{G_1, G_2\}$; $i = 1, \ldots, p$; $j = 1, \ldots, q$), a corresponding two-dimensional number $A_2$ can be written as:

$$A_2 = \alpha_{ij} \rfloor = \begin{matrix} \alpha_{pq} & \ldots & \alpha_{2q} & \alpha_{1q} \\ \ldots & \ldots & \ldots & \ldots \\ \alpha_{p2} & \ldots & \alpha_{22} & \alpha_{12} \\ \alpha_{p1} & \ldots & \alpha_{21} & \alpha_{11} \end{matrix} \rfloor.$$

By multisuccessor function it can be expressed as follows.

$$A_2 = \sum_{i,j} Sc_{ij}^{\alpha_{ij}}(\emptyset).$$

**Example 4.6.** Let us consider the polymset $A_2$, which has two-component basis of properties $G = \{G_1, G_2\}$:

A body form $G_1 = \{cube, pyramid, sphere, cone\}$;

A color $G_2 = \{black, white, red, green, blue\}$.

Suppose the polymset $A_2$ consist of one blue cube, eleven red pyramids, three black spheres, seven green cones, and five black cones. Then,

$A_2 = <(I_{15}, 1), (I_{23}, 11), (I_{31}, 3), (I_{44}, 7), (I_{41}, 5)> = <1_{15}, 11_{23}, 3_{31}, 7_{44}, 5_{41}>$.

$$A_2 = \begin{matrix} 0 & 0 & 0 & 1 \\ 7 & 0 & 0 & 0 \\ 0 & 0 & 11 & 0 \\ 0 & 0 & 0 & 0 \\ 5 & 3 & 0 & 0 \end{matrix} \rfloor. \quad \begin{matrix} \uparrow G_2 \\ \leftarrow G_1 \end{matrix}$$

A polymset $B_2$ be the immediate 2, 2-successor of the polymset $A_2$ if it consists of all $A_2$'s elements plus one white pyramid.

$B_2 = <(I_{15}, 1), (I_{23}, 11), (I_{31}, 3), (I_{44}, 7), (I_{41}, 5), (I_{22}, 1)> = <1_{15}, 11_{23}, 3_{31}, 7_{44}, 5_{41}, 1_{22}>$.

$$B_2 = \begin{matrix} 0 & 0 & 0 & 1 \\ 7 & 0 & 0 & 0 \\ 0 & 0 & 11 & 0 \\ 0 & 0 & 1 & 0 \\ 5 & 3 & 0 & 0 \end{matrix} \rfloor.$$

A polymset $D_2$ be the immediate 2, 3-predecessor of the polymset $A_2$ if it consists of all $A_2$'s elements without one red pyramid.

$D_2 = <(I_{15}, 1), (I_{23}, 10), (I_{31}, 3), (I_{44}, 7), (I_{41}, 5)> = <1_{15}, 10_{23}, 3_{31}, 7_{44}, 5_{41}>$.

$$D_2 = \begin{array}{cccc} 0 & 0 & 0 & 1 \\ 7 & 0 & 0 & 0 \\ 0 & 0 & 10 & 0 \\ 0 & 0 & 0 & 0 \\ 5 & 3 & 0 & 0 \end{array} \rfloor .$$

## 5. Peano system of natural multidimensional numbers

Let us introduce a concept of natural multidimensional numbers.
Here we consider two-dimensional numbers (2-numbers). It is also possible to develop a similar theory for multidimensional numbers with the higher than two dimensions.

**Definition 5.1.** The *natural 2-numbers* are called subpolymsets of a (universal) polymset $N_2$ such that for some polymsets $A_2$ and $B_2$ ($A_2, B_2 \in N_2$), there exists a multisuccessor function $B_2 = Sc_{ij}(A_2)$ and the following axioms hold.

I. There exists 2-number $\mathbf{0_2} = \mathbf{0_{ij}}\rfloor$, $i, j = 0, 1, 2, \ldots$, which is not a multisuccessor of any 2-numbers, i.e., $Sc_{pq}(A_2) \neq \mathbf{0_{ij}}\rfloor$, $\forall\ i, j, p, q$:

$$\mathbf{0_2} = \mathbf{0_{ij}}\rfloor = \begin{array}{cccc} 0_{pq} & \ldots & 0_{1q} & 0_{0q} \\ \ldots & \ldots & \ldots & \ldots \\ 0_{p1} & \ldots & 0_{11} & 0_{01} \\ 0_{p0} & \ldots & 0_{10} & 0_{00} \end{array} \rfloor$$

II. For any 2-number $A_2$ there exists the *unique* ij-successor $B_2 = Sc_{ij}(A_2)$, i.e., $A_2 = B_2 \Rightarrow Sc_{ij}(A_2) = Sc_{ij}(B_2)$. Meanwhile, $A_2 = B_2 \Rightarrow Sc_{ij}(A_2) \neq Sc_{pq}(B_2)$ if $i \neq p$, $j \neq q$. It follows that the multisuccessor function is one-to-one correspondence.

III. Any S-number is an ij-successor for no more than one S-number, i.e.,
 $Sc_{ij}(A_2) = Sc_{ij}(B_2) \Rightarrow A_2 = B_2$.

IV. /The induction axiom/. Any set $\mathbf{P_2}$ of 2-numbers contains all natural 2-numbers, i.e., $\mathbf{P_2}$ is equal to $\mathbf{N_2}$, if the following conditions hold:
 (i) $\mathbf{0_{ij}}\rfloor \in \mathbf{P_2}$,
 (ii) $A_2 \in \mathbf{P_2} \Rightarrow Sc_{ij}(A_2) \in \mathbf{P_2}$.

The axioms I-IV are a generalization of well-known Peano axioms, so the Definition 5.1 allows to define Peano 2-system ($\mathbf{P_2}$, $Sc_{ij}$, $\mathbf{0_{ij}}\rfloor$) of natural 2-numbers and to suggest the axiom of existence.

**Axiom 5.2.** There exists at least one Peano 2-system.

**Theorem 5.3.** Any two Peano 2-systems are isomorphic.

**Theorem 5.4.** Let ($\mathbf{P_2}$, $Sc_{ij}$, $\mathbf{0_{ij}}\rfloor$) be Peano 2-system. Then for any $A_2 \in \mathbf{P_2}$ either $A_2 = \mathbf{0_{ij}}\rfloor$ or there exists such $B_2 \in \mathbf{P_2}$ that $A_2 = Sc_{ij}(B_2)$, and this $B_2$ is unique for given $i, j$.

**Proof.** To prove this theorem we use the induction axiom (IV).

Let $\mathbf{A_2} = \{A_2: A_2 \in \mathbf{P_2}\}$, and $A_2 = 0_{ij}\rfloor$ or $A_2 = Sc_{ij}(B_2)$ for some $B_2 \in \mathbf{P_2}$. Prove that $\mathbf{A_2} = \mathbf{P_2}$.

It is clear that $\mathbf{A_2} \subseteq \mathbf{P_2}$ and $0_{ij}\rfloor \in \mathbf{A_2}$. Show that $Sc_{ij}(A_2) \in \mathbf{A_2}, \forall A_2 \in \mathbf{A_2}$, i.e. for any $A_2 \in \mathbf{A_2}$ either $Sc_{ij}(A_2) = 0_{ij}\rfloor$ (that is impossible by axiom I) or $Sc_{ij}(A_2) = Sc_{ij}(C_2)$ for some $C_2 \in \mathbf{P_2}$ (it may be true supposing $A_2 = C_2$). Then $\mathbf{A_2} = \mathbf{P_2}$.

Uniqueness follows from axiom III. From $Sc_{ij}(B_2) = A_2$ and $Sc_{ij}(C_2) = A_2$ it follows $B_2 = C_2$.

Hereinafter suppose that $(\mathbf{P_2}, Sc_{ij}, 0_{ij}\rfloor)$ is fixed Peano 2-system, and $\mathbf{P_2}$ is called a *set of natural 2-numbers*.

By using axioms I-IV it is easy to prove the next useful theorems.

**Theorem 5.5.** If ij-successors (2-numbers) of some initial 2-numbers are distinct, then the initial 2-numbers are distinct: $Sc_{ij}(A_2) \neq Sc_{ij}(B_2) \Rightarrow A_2 \neq B_2$.

Note that for distinct index pairs this theorem is not true, i.e., not always from $Sc_{ij}(A_2) \neq Sc_{pq}(B_2)$ it follows $A_2 \neq B_2$.
Consider $A_2 = B_2$, then $Sc_{ij}(A_2) = C_2$, $Sc_{pq}(A_2) = D_2$. By axiom II, $D_2 \neq C_2$.

**Theorem 5.6.** If given 2-numbers are distinct, then their ij-successors are different: $A_2 \neq B_2 \Rightarrow Sc_{ij}(A_2) \neq Sc_{ij}(B_2)$.

Note that for distinct index pairs this theorem is not true, i.e., for i≠p, j≠q it is possible to obtain $A_2 \neq B_2 \Rightarrow Sc_{ij}(A_2) = Sc_{pq}(B_2)$.

**Example 5.7.**  
Let $A_2 = 10\overset{00}{\rfloor}$, $B_2 = 00\overset{01}{\rfloor}$, i.e., $A_2 \neq B_2$.

Then $Sc_{01}(A_2) = 10\overset{01}{\rfloor}$, and $Sc_{10}(B_2) = 10\overset{01}{\rfloor}$, i.e., $Sc_{01}(A_2) = Sc_{10}(B_2)$.

**Theorem 5.8.** Any 2-number is distinct from a 2-number successeding it, i.e., $Sc_{ij}(A_2) \neq A_2, \forall i, j, \forall A_2 \in \mathbf{P_2}$.

**Proof.** Let $\mathbf{M}$ be a set of natural 2-numbers for which the theorem is true.
1). By axiom I, $Sc_{pq}(1_{ij}) \neq 1_{ij}$, so $1_{ij} \in \mathbf{M}$;
2). If $A_2 \in \mathbf{M}$, then $Sc_{ij}(A_2) \neq A_2$. Hence (by theorem 5.4.) $Sc_{ij}(Sc_{pq}(A_2)) \neq Sc_{pq}(A_2)$, i.e., $Sc_{pq}(A_2) \in \mathbf{M}$.

By axiom IV, $\mathbf{M}$ contains all 2-numbers, i.e., $A_2 \neq Sc_{ij}(A_2), \forall A_2$.

## 6. Addition of natural 2-numbers

In this paper, following [4], we introduce main definitions and prove basic theorems.

**Definition 6.1.** The *addition* of natural 2-numbers is defined as a correspondence which for any pair of 2-numbers $A_2$ and $B_2$ is compared with an unique 2-number $A_2 + B_2$ and the following conditions hold:

(i)  + is a binary operation on $\mathbf{P_2}$ under which $\mathbf{P_2}$ is closed;

(ii) $A_2 + 1_{ij} = Sc_{ij}(A_2)$, $\forall\, A_2 \in \mathbf{P_2}$;
(iii) $A_2 + Sc_{ij}(B_2) = Sc_{ij}(A_2 + B_2)$, $\forall\, A_2, B_2 \in \mathbf{P_2}$.

**Theorem 6.2.** There *exists* the addition of natural 2-numbers and it is *unique*. That is we have one and only one correspondence which is compared with any pair of 2-numbers $A_2$ and $B_2$ the 2-numbers $A_2 + B_2$ such that:

1. $A_2 + 1_{ij} = Sc_{ij}(A_2)$, $\forall\, A_2 \in \mathbf{P_2}$;
2. $A_2 + Sc_{ij}(B_2) = Sc_{ij}(A_2 + B_2)$, $\forall\, A_2, B_2 \in \mathbf{P_2}$.

**Proof.** See [3].

**Theorem 6.3.** (*Associative law for* +) For any natural 2-numbers $A_2, B_2, C_2 \in \mathbf{P_2}$, we have
$(A_2 + B_2) + C_2 = A_2 + (B_2 + C_2)$.

**Proof.** By induction on $C_2$. Consider any $A_2, B_2 \in \mathbf{P_2}$. Let $\mathbf{M}$ ($\mathbf{M} \subseteq \mathbf{P_2}$) be a set of 2-numbers $C_2$ for which the equality is true. Then

(i) $(A_2 + B_2) + 1_{ij} = Sc_{ij}(A_2 + B_2) = A_2 + Sc_{ij}(B_2) = A_2 + (B_2 + 1_{ij})$, $1_{ij} \in \mathbf{M}$.
(ii) If $C_2 \in \mathbf{M}$ then $(A_2 + B_2) + C_2 = A_2 + (B_2 + C_2)$.
Hence $(A_2 + B_2) + Sc_{ij}(C_2) = Sc_{ij}((A_2 + B_2) + C_2) = A_2 + Sc_{ij}(B_2 + C_2) = A_2 + (B_2 + Sc_{ij}(C_2))$, that is $Sc_{ij}(C_2) \in \mathbf{M}$. Therefore by axiom IV the equality $(A_2 + B_2) + C_2 = A_2 + (B_2 + C_2)$ is true for any $A_2, B_2, C_2 \in \mathbf{P_2}$.

**Theorem 6.4.** (*Commutative law for* +) For any natural 2-numbers $A_2, B_2 \in \mathbf{P_2}$, it follows that
$A_2 + B_2 = B_2 + A_2$.

**Proof**. Let us show by induction on $A_2$ that $A_2 + 1_{ij} = 1_{ij} + A_2$, $\forall\, A_2 \in \mathbf{P_2}$.
Let $\mathbf{M}$ ($\mathbf{M} \subseteq \mathbf{P_2}$) be a set of 2-numbers $A_2$ for which this equality is true. It is clear that $1_{ij} \in \mathbf{M}$. Moreover, if $A_2 \in \mathbf{M}$, then $A_2 + 1_{ij} \in \mathbf{M}$.
Hence $Sc_{ij}(A_2) + 1_{ij} = (A_2 + 1_{ij}) + 1_{ij} = (1_{ij} + A_2) + 1_{ij} = Sc_{ij}(1_{ij} + A_2) = 1_{ij} + Sc_{ij}(A_2)$, i.e., $Sc_{ij}(A_2) \in \mathbf{M}$. By axiom IV we get $A_2 + 1_{ij} = 1_{ij} + A_2$.
Let us show by induction on $B_2$ that $A_2 + B_2 = B_2 + A_2$. Let $\mathbf{M}$ ($\mathbf{M} \subseteq \mathbf{P_2}$) be a set of 2-numbers $B_2$ for which this equality is true for given $A_2$. Given $1_{ij} \in \mathbf{M}$, if $B_2 \in \mathbf{M}$, then $A_2 + B_2 = B_2 + A_2$. By using theorem 5.4 we obtain
$A_2 + Sc_{ij}(B_2) = Sc_{ij}(A_2 + B_2) = Sc_{ij}(B_2 + A_2) = B_2 + Sc_{ij}(A_2) = B_2 + (A_2 + 1_{ij}) =$
$= B_2 + (1_{ij} + A_2) = (B_2 + 1_{ij}) + A_2 = Sc_{ij}(B_2) + A_2$, i.e., $Sc_{ij}(B_2) \in \mathbf{M}$.

**Theorem 6.5.** For any natural 2-numbers $A_2, B_2 \in \mathbf{P_2}$ ($B_2 \neq \mathbf{0_2}$), $A_2 + B_2 \neq B_2$.

**Proof.** This theorem is true for $B_2 = 1_{ij}$ because $A_2 + 1_{ij} = Sc_{ij}(A_2) \neq 1_{ij}$ by axiom I. If $A_2 + B_2 \neq B_2$, then by theorem 5.6: $A_2 + Sc_{ij}(B_2) = Sc_{ij}(A_2 + B_2) \neq Sc_{ij}(B_2)$.

We can prove a cancellation law for + by analogy.

**Theorem 6.6.** (*Cancellation law for* + ) For any $A_2, B_2, C_2 \in \mathbf{P_2}$ from $A_2 + C_2 = B_2 + C_2$ it follows that $A_2 = B_2$.

**Theorem 6.7**. (*Tetratomy* law for +) For any natural 2-numbers $A_2, B_2 \in \mathbf{P_2}$, one of the following four cases holds:

(i) $\quad A_2 = B_2$.
(ii) $\quad$ For some $C_2 \in P_2$, $A_2 = B_2 + C_2$.
(iii) $\quad$ For some $D_2 \in P_2$, что $B_2 = A_2 + D_2$.
(iv) $\quad A_2$ and $B_2$ are not comparable ($A_2 \lozenge B_2$), i.e., $A_2$ and $B_2$ cannot be received one of another by means of the addition operation.

**Proof**. By theorem 6.5 cases (i) and (ii), (i) and (iii) cannot hold simultaneously.
If case (ii) holds with case (iii), then $A_2 = B_2 + C_2 = (A_2 + D_2) + C_2 = A_2 + (D_2 + C_2) \neq A_2$ that contradicts the theorem 6.5.

Cases (i) and (iv), (ii) and (iv), (iii) and (iv) are not compatible by the definition.
Prove that one of cases (i)-(iv) always takes place.
Let $\mathbf{M}$ ($\mathbf{M} \subseteq \mathbf{P_2}$) be a set of 2-numbers $B_2$ such that for given $A_2$ one of cases (i)-(iv) takes place.

(a) If $A_2 = 1_{ij}$, then case (i) holds for $B_2 = 1_{ij}$. If $A_2 \neq 1_{ij}$, then
$A_2 = Sc_{ij}(Q_2) = Q_2 + 1_{ij}$, i.e., case (ii) holds for $B_2 = 1_{ij}$ and $Q_2 \neq 0$ such that $1_{ij} \in \mathbf{M}$;

(b) Let $B_2 \in \mathbf{M}$. Then either $A_2 = B_2$ that is $Sc_{ij}(B_2) = B_2 + 1_{ij}$ (case (iii) for $Sc_{ij}(B_2)$),
or $A_2 = B_2 + K_2$, and, if $K_2 = 1_{ij}$, then $A_2 = B_2 + 1_{ij} = Sc_{ij}(B_2)$ (case (i) for $Sc_{ij}(B_2)$),
but if $K_2 \neq 1_{ij}$, then $K_2 = Sc_{pq}(L_2)$, and $A_2 = B_2 + Sc_{pq}(L_2) = B_2 + (L_2 + 1_{pq}) = B_2 + (1_{pq} + L_2) = (B_2 + 1_{pq}) + L_2 = Sc_{pq}(B_2) + L_2$, i.e., case (ii) for $Sc_{ij}(B_2)$,
or $B_2 = A_2 + L_2$, and $Sc_{ij}(B_2) = Sc_{ij}(A_2 + L_2) = A_2 + Sc_{ij}(L_2)$, i.e., case (iii) for $Sc_{ij}(B_2)$.
Hence $B_2 \in \mathbf{M}$;

(c) For the given 2-number $A_2$ let us show that there exists 2-number $B_2 \in \mathbf{M}$ such that $A_2 \neq B_2$, $A_2 \neq B_2 + C_2$, $B_2 \neq A_2 + D_2$, $\forall C_2, D_2 \in \mathbf{M}$.
Let there exist 2-number $G_2$ such that case (ii) holds both for $A_2$ and for $B_2$. $G_2 = A_2 + H_2$, $G_2 = B_2 + J_2$.
Then $A_2 + H_2 = B_2 + J_2$, $A_2 = B_2 + J_2 - H_2$ or $A_2 = B_2 + (J_2 - H_2)$. But the subtraction operation was not defined on the set $\mathbf{N_2}$.
Hence $A_2 \neq B_2$ и $A_2 \neq B_2 + C_2$.
Meanwhile, from $B_2 + J_2 = A_2 + H_2$ follows $B_2 = A_2 + H_2 - J_2$, $B_2 = A_2 + (H_2 - J_2)$. Finally, we obtain $B_2 \neq A_2$ and $B_2 \neq A_2 + D_2$.

By multisuccessor function we can assert that if $A_2 = \sum_{i,j} Sc_{ij}^{\alpha_{ij}}(\varnothing)$ and $B_2 = \sum_{i,j} Sc_{ij}^{\beta_{ij}}(\varnothing)$, then $A_2 + B_2 = \sum_{i,j} Sc_{ij}^{\alpha_{ij}+\beta_{ij}}(\varnothing)$.

Argumentation mentioned above allows to define a new algebraic structure (commutative semigroup (monoid)), i.e., the set of natural 2-numbers with binary operation +, neutral element (identity) $\mathbf{0_2} := \mathbf{0_{ij}}\rfloor$ and two-side cancellation.

## 7. Multiplication of natural 2-numbers

Let us introduce multiplication of natural 2-numbers.

**Proposition 7.1.** There exists a binary operation · on the set of natural 2-numbers such that

$$Sc_{ij}(\emptyset) \cdot Sc_{pq}(\emptyset) = Sc_{i+p\ j+q}(\emptyset).$$

In 2-number representation we obtain $1_{ij} \cdot 1_{pq} = 1_{i+p\ j+q}$.

**Proposition 7.2.** If $A_2 = \sum_{i,j} Sc_{ij}^{\alpha_{ij}}(\emptyset)$ and $B_2 = \sum_{i,j} Sc_{pq}^{\beta_{pq}}(\emptyset)$, then

$$A_2 \cdot B_2 = \sum_{i,j} Sc_{i+p\_j+q}^{\alpha_{ij} \cdot \beta_{pq}}(\emptyset).$$

In 2-number representation we get $\alpha_{ij} \cdot \beta_{pq} = (\alpha \cdot \beta)_{i+p\ j+q}$.

**Definition 7.3.** The *multiplication* of natural 2-numbers is defined as a correspondence that for any pair of 2-numbers $A_2$ and $B_2$ assigns a unique 2-number $A_2 \cdot B_2$, and the following conditions hold:

(i)   · is a binary operation on $\mathbf{P_2}$ under which $\mathbf{P_2}$ is closed;
(ii)  $A_2 \cdot 1_{ij} = (A_2)_{ij},\ \forall\ A_2 \in \mathbf{P_2}$;
(iii) $A_2 \cdot Sc_{ij}(B_2) = A_2 \cdot B_2 + (A_2)_{ij},\ \forall\ A_2, B_2 \in \mathbf{P_2}$.

By $(.)_{ij}$ we denote a shift operation. This operation shifts a given 2-number by $i$ positions to the left and $j$ positions upwards.

**Example 7.4.**

Let $A_2 = \begin{matrix} \phantom{0}02 \\ 21 \rfloor \end{matrix}$, $(.)_{ij} = (.)_{11}$, then $(A_2)_{11} = \begin{matrix} 020 \\ 210 \\ 000 \rfloor \end{matrix}$.

**Theorem 7.5.** There *exists* the multiplication of natural 2-numbers and it is *unique*.

**Proof.** See [3].

**Theorem 7.6.** (*Right distributive law for* · *over* +) For any $A_2, B_2, C_2 \in \mathbf{P_2}$,

$$(A_2 + B_2)\, C_2 = A_2 \cdot C_2 + B_2 \cdot C_2.$$

**Proof** is performed by induction on $C_2$.
  (i) $(A_2 + B_2)\, 1_{ij} = (A_2 + B_2)_{ij} = (A_2)_{ij} + (B_2)_{ij} = A_2 \cdot 1_{ij} + B_2 \cdot 1_{ij}$, that is for $C_2 = 1_{ij}$ the theorem is true.
  (ii) If the theorem is true for $C_2$, then $(A_2 + B_2)\, C_2 = A_2 \cdot C_2 + B_2 \cdot C_2$. By use the associative law and the commutative law for +, we get for $Sc_{ij}(C_2)$:

$$(A_2 + B_2)\, Sc_{ij}(C_2) = (A_2 + B_2)\, C_2 + (A_2 + B_2)_{ij} =$$
$$= (A_2\, C_2 + B_2\, C_2) + ((A_2)_{ij} + (B_2)_{ij}) = (A_2 \cdot C_2 + (A_2)_{ij}) + (B_2 \cdot C_2 + (B_2)_{ij}) =$$

$$= (A_2 \cdot C_2 + A_2 \cdot 1_{ij}) + (B_2 \cdot C_2 + B_2 \cdot 1_{ij}) = A_2 (C_2 + 1_{ij}) + B_2 (C_2 + 1_{ij}) =$$
$$= A_2 \cdot Sc_{ij}(C_2) + B_2 \cdot Sc_{ij}(C_2).$$

**Theorem 7.7.** (*Commutative law for* ·) For any $A_2, B_2, C_2 \in \mathbf{P_2}$,
$$A_2 \cdot B_2 = B_2 \cdot A_2.$$

**Proof**. (i) By induction on $B_2$, we prove this theorem when $A_2 = 1_{ij}$, i.e., $1_{ij} \cdot B_2 = B_2 \cdot 1_{ij}$.

Let **M** be a set of 2-numbers $B_2$ for which the equality from the theorem is true. Then $1_{ij} \in \mathbf{M}$.

If $1_{ij} \cdot B_2 = B_2 \cdot 1_{ij}$, then $1_{ij} \cdot Sc_{pq}(B_2) = 1_{ij} \cdot B_2 + 1_{ij} \cdot 1_{pq} = B_2 \cdot 1_{ij} + 1_{pq} \cdot 1_{ij} = (B_2)_{ij} + (1_{pq})_{ij} = (Sc_{pq}(B_2))_{ij} = Sc_{pq}(B_2) \cdot 1_{ij}$, i.e., $Sc_{pq}(B_2) \in \mathbf{M}$.

(ii) By induction on $A_2$, we prove that $A_2 \cdot B_2 = B_2 \cdot A_2$ for the given $B_2$ and $A_2 \in \mathbf{M}$, where $A_2 \cdot B_2 = B_2 \cdot A_2$.

By (i) $1_{ij} \in \mathbf{M}$. If $A_2 \in \mathbf{M}$, then $A_2 \cdot B_2 = B_2 \cdot A_2$ and by theorem 7.4:
$$Sc_{ij}(A_2) \cdot B_2 = (A_2 + 1_{ij}) B_2 = A_2 \cdot B_2 + 1_{ij} \cdot B_2 = B_2 \cdot A_2 + B_2 \cdot 1_{ij} =$$
$$= B_2 (A_2 + 1_{ij}) = B_2 \cdot Sc_{ij}(A_2), \text{ hence } Sc_{ij}(A_2) \in \mathbf{M}.$$

Using theorems 7.4, 7.5, it is easy to prove:

**Theorem 7.8**. (*Left distributive law*) For any $A_2, B_2, C_2 \in \mathbf{P_2}$,
$$C_2 (A_2 + B_2) = C_2 \cdot A_2 + C_2 \cdot B_2.$$

**Theorem 7.9**. (*Associative law for* ·) For any $A_2, B_2, C_2 \in \mathbf{P_2}$,
$$(A_2 \cdot B_2) C_2 = A_2 (B_2 \cdot C_2).$$

**Proof.** Consider any $A_2, B_2 \in \mathbf{P_2}$ and let **M** be a set of 2-numbers for which the conditions are true.

(i) $(A_2 \cdot B_2) 1_{ij} = A_2 (B_2 \cdot 1_{ij})$, $1_{ij} \in \mathbf{M}$.
(ii) If $C_2 \in \mathbf{M}$, then $(A_2 \cdot B_2) C_2 = A_2 (B_2 \cdot C_2)$.

Using theorem 7.6, we obtain $(A_2 \cdot B_2) Sc_{ij}(C_2) = (A_2 \cdot B_2) C_2 + (A_2 \cdot B_2) 1_{ij} =$
$$= A_2 (B_2 \cdot C_2) + A_2 (B_2 \cdot 1_{ij}) = A_2 (B_2 \cdot C_2 + B_2 \cdot 1_{ij}) = A_2 (B_2 (C_2 + 1_{ij})) =$$
$$= A_2 (B_2 \cdot Sc_{ij}(C_2)), Sc_{ij}(C_2) \in \mathbf{M}.$$

We prove that $1_{00}$ is the absolute unity for the Peano 2-system of natural 2-numbers $\mathbf{P_2}$.

**Lemma 7.10.** $1_{00} \cdot A_2 = A_2$.

**Proof.** By induction on $A_2$. $1_{00} \cdot 1_{ij} = 1_{i+0\,j+0} = 1_{ij}$.
Let $1_{00} \cdot A_2 = A_2$, then $1_{00} \cdot Sc_{ij}(A_2) = (1_{00} \cdot A_2) + 1_{ij} = A_2 + 1_{ij} = Sc_{ij}(A_2)$.

The following lemma is similarly proved.

**Lemma 7.11.** $1_{ij} \cdot A_2 = (A_2)_{ij}$, $\forall A_2 \in \mathbf{P_2}$.

Thus, the set of natural 2-numbers equipped with two binary operations + and · is a commutative semiring $\langle \mathbf{P_2}, +, \cdot, \mathbf{0_2}, 1_{00} \rangle$.

## 8. Conclusion

In this paper, we have proposed the concept of natural multidimensional numbers, many related notions have been defined, and many related theorems have been proved.

Future studies of the new algebraic structure may lead to the m-numbers ordering determination and extension of the natural multidimensional numbers semiring to the ring of whole multidimensional numbers.

Moreover, the concept of multidimensional numbers is important and useful for the definition and development of special multidimensional numeration systems.

## 9. Acknowledgement

The author is grateful to Dr. Mark Burgin for useful remarks and suggestions.

## References

[1] R. Bruni, F. Gadducci, Some algebraic laws for spans (and their connections with multirelations), Electronic Notes in Theoretical Computer Science 44 (3) (2003) 175-193.

[2] C. Calude et al., Multiset Processing, Lecture Notes in Computer Science 235 (2001).

[3] A. Chunikhin, Introduction to Multidimensional Number Systems. Theoretical Foundations and Applications, Kiev, 2006 (in Russian).

[4] S. Feferman, The Number Systems. Foundations of Algebra and Analysis, Addison-Wesley, Reading, 1964.

[5] L. Feijs, R.L. Krikhaar, Relation algebra with multi-relations, International Journal of Computer Mathematics 70 (1998) 57–74.

[6] S.C. Kleene, Introduction to Metamathematics, D. Van Nostrand Co., Inc., 1952.

[7] A.J. Knobbe, Multi-Relational Data Mining, IOS Press, 2006.